\documentclass[final]{siamltex}

\usepackage{amsmath} 
\usepackage{amsfonts}

\usepackage{multirow}

\title{An Augmented Lagrangian Preconditioner for Navier--Stokes Equations with Runge--Kutta in Time}

\author{Santolo Leveque\footnotemark[1] \and Yunhui He\footnotemark[1] \and Maxim Olshanskii\footnotemark[1]}

\begin{document}
\maketitle

\footnotetext[1]{Department of Mathematics, University of Houston, 3551 Cullen Blvd, Houston, Texas 77204-3008, USA ({\tt sleveque@central.uh.edu, yhe43@central.uh.edu, maolshanskiy@uh.edu})}

\begin{abstract}
We consider a Runge--Kutta method for the numerical time integration of the nonstationary incompressible Navier--Stokes equations. This yields a sequence of nonlinear problems to be solved for the stages of the Runge--Kutta method. The resulting nonlinear system of differential equations is discretized using a finite element method. To compute a numerical approximation of the stages at each time step, we employ Newton's method, which requires the solution of a large and sparse generalized saddle-point problem at each nonlinear iteration. We devise an augmented Lagrangian preconditioner within the flexible GMRES method for solving the Newton systems at each time step. The preconditioner can be applied inexactly with the help of a multigrid routine. We present numerical evidence of the robustness and efficiency of the proposed strategy for different values of the viscosity, mesh size, time step, and number of stages of the Runge--Kutta method.
\end{abstract}



\pagestyle{myheadings}
\thispagestyle{plain}
\markboth{S. LEVEQUE, Y. HE, AND M. OLSHANSKII}{AL PRECONDITIONER FOR NAVIER--STOKES RUNGE--KUTTA IN TIME}

\section{Introduction}\label{sec_1}
Accurately modeling the motion of Newtonian incompressible viscous fluids is a fundamental problem  in computational science and engineering. The governing equations for such flows are the incompressible Navier–Stokes equations, whose numerical solution has been the focus of extensive research for decades. Many time-stepping schemes for the nonstationary Navier–Stokes equations require solving a nonlinear problem at each time step, typically via Picard or Newton iterations. These nonlinear solvers in turn generate a sequence of linear systems that must be solved efficiently at each iteration.

Classical time integration schemes for the nonstationary incompressible Navier--Stokes equations include multistep methods such as backward Euler and the trapezoidal rule \cite{Temam, William}, projection methods \cite{Guermond_Minev_Shen}, and implicit-explicit (IMEX) methods \cite{Marion_Roger}, among others. There are many good reasons for using these methods. First, the block structure of the discretized system allows one to adapt linear solvers for the stationary Navier--Stokes problem to the nonstationary case. Furthermore, some of these methods are A-stable, so no Courant--Friedrichs--Lewy (CFL) condition must be satisfied, meaning that no restriction on the time step is required. These methods satisfy a proper energy balance that ensures energy stability.

Despite the advantages described above, classical time integrators for the nonstationary incompressible Navier--Stokes equations are typically of order no greater than two. Specifically, for linear multistep methods, the second Dahlquist barrier states that an A-stable linear multistep method cannot have an order of convergence greater than two; see, for example, \cite[Theorem 6.6]{Lambert}. However, the work in \cite{Huang_Shen} develops high-order (up to fifth-order) A-stable IMEX methods for the nonstationary incompressible Navier--Stokes equations under periodic boundary conditions.

In contrast to linear multistep methods, it is possible to devise A-stable implicit Runge--Kutta methods of arbitrary order. Moreover, some implicit Runge--Kutta methods can be constructed to be L- or B-stable; see, for instance, \cite{Butcher, Hairer_Wanner, Lambert}. This motivates us to consider Runge--Kutta methods for the nonstationary incompressible Navier--Stokes equations. Although there are several studies applying Runge--Kutta methods to the time integration of the Navier--Stokes equations, we are aware of only \cite{AbuLabdeh_MachLachlan_Farrell} that addresses the linear solver aspect of this approach. Specifically, the authors of \cite{AbuLabdeh_MachLachlan_Farrell} designed a monolithic multigrid method as the stage solver.

The recent rediscovery of Runge--Kutta methods as time integrators for the numerical solution of PDEs has led to a plethora of linear solvers for the system of stages of the method. After the development of simple iterations such as the block-Jacobi \cite{Mardal_Nilssen_Staff} and block-Gauss--Seidel \cite{Staff_Mardal_Nilssen} preconditioners for the heat equation, more advanced techniques have been employed to find good approximations of the stage systems. For example, strategies based on special LU factorizations of the coefficient matrix \cite{Rana_Howle_Long_Meek_Milestone} and diagonalizations of its lower triangular factor \cite{Axelsson_Dravins_Neytcheva} have been proposed. Further, in \cite{Munch_Dravins_Kronbichler_Neytcheva}, the authors develop a stage-parallel implementation of the preconditioner derived in \cite{Axelsson_Dravins_Neytcheva}. Regarding parallel implementations of Runge--Kutta integrators for PDEs, we mention \cite{Kressner_Massei_Zhu, Leveque_Bergamaschi_Martinez_Pearson}, where the authors propose parallel-in-time solvers for the all-at-once space-time discretization of the heat equation with Runge--Kutta in time. The solver in \cite{Leveque_Bergamaschi_Martinez_Pearson} is further adapted to handle the incompressible Stokes equations. Other examples of stage solvers include methods based on the eigenvalue decomposition of the inverse of the coefficient matrix \cite{Southworth_Krzysik_Pazner_Sterck_1}, the Schur decomposition of the coefficient matrix \cite{Southworth_Krzysik_Pazner_Sterck_2}, and monolithic multigrid solvers \cite{AbuLabdeh_MachLachlan_Farrell, Farrell_Kirby_MarchenaMenendez, Kirby}. We also mention early works on monolithic multigrid methods for stage solvers by Vandewalle and co-authors \cite{Boonen_VanLent_Vanderwalle, Rosseel_Boonen_Vanderwalle, VanLent_Vanderwalle}.

The preconditioning strategy we adopt here is based on an augmented Lagrangian approach. Originally, the augmented Lagrangian preconditioner was introduced in \cite{Benzi_Olshanskii_2006} for solving the Picard linearization of the stationary incompressible Navier--Stokes equations. This preconditioning strategy was shown to be robust with respect to the mesh size and the viscosity of the fluid---features that are particularly important for solver performance. The key components of the preconditioner are a specialized multigrid routine used for the approximate inversion of the augmented momentum equation and a special approximation of the Schur complement, derived from the application of the Sherman--Morrison--Woodbury formula. In \cite{Benzi_Olshanskii_2011}, the authors proved the optimality of the preconditioner via a field-of-values analysis of the preconditioned Schur complement, showing that its eigenvalues are contained in a rectangle in the right half of the complex plane, with boundaries independent of the viscosity, provided the augmentation parameter is proportional to the inverse of the viscosity. However, as demonstrated by the numerical results in \cite{Benzi_Olshanskii_2006}, robust convergence is also achieved when the augmentation parameter is of order one, provided a suitable rescaling of the equations is applied. In \cite{Farrell_Mitchell_Wechsung}, the authors apply this augmented Lagrangian strategy to the Newton linearization of the 3D stationary incompressible Navier--Stokes equations. In this case as well, the numerical results confirm the effectiveness and robustness of the approach. Finally, we mention that in \cite{Leveque_Benzi_Farrell}, an augmented Lagrangian-based preconditioner is employed in an inexact Newton iteration for the stationary incompressible Navier--Stokes control problem.

As we will show in the numerical results, the augmented Lagrangian-based preconditioner we adopt in this work is robust with respect to both the problem parameters and the choice of Runge--Kutta method, resulting in a linear solver that scales linearly with the problem size.

This paper is organized as follows. In Section~\ref{sec_2}, we introduce Runge--Kutta methods, which are then applied to the nonstationary incompressible Navier--Stokes equations in Section~\ref{sec_3}. In Section~\ref{sec_4}, we present the augmented Lagrangian preconditioner used to solve the stage systems of the Runge--Kutta method. Numerical evidence of the efficiency and robustness of the preconditioner is given in Section~\ref{sec_5}. Finally, conclusions and directions for future work are discussed in Section~\ref{sec_6}.

\section{Runge--Kutta Methods}\label{sec_2}
We recall a variant of the Runge--Kutta method for an index 2 differential--algebraic equation (DAE) of the form  
\begin{equation}\label{general_DAE}
\left\{
\begin{array}{ll}
\vspace{0.25ex}
\dot{u} = \mathcal{F}(\hat{w},t) & \quad \mathrm{in} \; (0,T), \\
\vspace{0.25ex}
\mathcal{G}(u,t) = 0 & \quad \mathrm{in} \; (0,T),
\end{array}
\right.
\end{equation}  
with $\hat{w} = (u,p)$ and a suitable initial condition $\hat{w}_0=(u_0,p_0)$.  

A Runge--Kutta method is fully characterized by the nodes $\mathbf{c}_{\mathrm{RK}} :=\left\{ c_i \right\}_{i=1}^s$, the weights $\mathbf{b}_{\mathrm{RK}}:=\left\{ b_i \right\}_{i=1}^s$ of the quadrature rule employed, and the coefficients $A_{\mathrm{RK}} :=\left\{ a_{i,j} \right\}_{i,j=1}^s$ of the method, and it can be represented in compact form by the following Butcher tableau:
\begin{displaymath}
\def\arraystretch{1.2}
\begin{array}{c|c}
\mathbf{c}_{\mathrm{RK}} & A_{\mathrm{RK}}\\
\hline
 & \mathbf{b}_{\mathrm{RK}}^\top
\end{array}
\end{displaymath}

After dividing the time interval $[0, T]$ into $n_t$ equally spaced subintervals of length $\Delta t$, a time-stepping technique based on an $s$-stage Runge--Kutta method updates the numerical solution of the ODE $\dot{u}(t)=\mathcal{F}(u(t),t)$ as follows:
\begin{displaymath}
u_{n+1}= u_{n} + \Delta t \sum_{i=1}^{s} b_i Y_{i,n}, \quad n = 0, \ldots, n_t-1.
\end{displaymath}
The quantities $Y_{i,n}$, for $i = 1, \ldots, s$, are called the \emph{stages} between $t_n:= n \Delta t$ and $t_{n} + \Delta t$, for $n=0, \ldots, n_t-1$, and solve the following nonlinear equations:
\begin{equation}\label{stages_Runge_Kutta}
Y_{i,n} = \mathcal{F}\left(u_n + \Delta t \sum_{j=1}^{s} a_{i,j} Y_{j,n}, t_n + c_i \Delta t\right), \qquad i=1,\ldots,s.
\end{equation}

In what follows, we assume that the matrix $A_{\mathrm{RK}}$ is invertible. This assumption is common in many works on preconditioners for the stage system and is not restrictive. For example, one can prove that for an implicit Runge--Kutta method with invertible $A_{\mathrm{RK}}$, if $a_{s,j}=b_j$ for $j=1,\ldots,s$, or $a_{i,1} = b_1$ for $i=1,\ldots,s$, then the method is L-stable; see \cite[Proposition 3.8, p. 45]{Hairer_Wanner}.

To apply a Runge--Kutta method in time to solve \eqref{general_DAE}, we introduce the auxiliary variables $\hat{w}_{i,n} = [w^{u}_{i,n}, w^{p}_{i,n}]$, where
\begin{equation}\label{wu_wp}
w^{u}_{i,n} = u_n + \Delta t \sum_{j=1}^{s} a_{i,j} Y^{u}_{j,n}, \qquad w^{p}_{i,n} = p_n + \Delta t \sum_{j=1}^{s} a_{i,j} Y^{p}_{j,n}.
\end{equation}
Then, for the index 2 DAE in \eqref{general_DAE} equation \eqref{stages_Runge_Kutta} yields to 
\begin{equation}\label{stages_Runge_Kutta_2}
\left\{
\begin{array}{ll}
\vspace{0.25ex}
Y^{u}_{i,n} = \mathcal{F}\left(\hat{w}_{i,n}, t_n + c_i \Delta t\right), & \quad i=1,\ldots,s,\\
\mathcal{G}(w^{u}_{i,n},t) = 0, & \quad i=1,\ldots,s.
\end{array}
\right.
\end{equation}

The nonlinearity of the problem requires solving a linearized version of \eqref{stages_Runge_Kutta_2} iteratively to obtain the stages $Y^{u}_{i,n}$ and $Y^{p}_{i,n}$, for $i=1, \ldots, s$. A natural approach is to apply Newton's method to \eqref{stages_Runge_Kutta_2}.

Once the solutions $\left\{ Y^{u}_{i,n} \right\}_{i=1}^s$ and $\left\{ Y^{p}_{i,n} \right\}_{i=1}^s$ of \eqref{stages_Runge_Kutta_2} are found, the numerical approximations at time $t_{n+1}$ are computed using the following recurrence formula:
\begin{equation}\label{DAE_RK_step}
u_{n+1}= u_{n} + \Delta t \sum_{i=1}^{s} b_i Y^{u}_{i,n}, \qquad p_{n+1}= p_{n} + \Delta t \sum_{i=1}^{s} b_i Y^{p}_{i,n}.
\end{equation}

We conclude this section with some remarks on the initial condition $\hat{w}_0$. As mentioned, the system \eqref{general_DAE} must be equipped with a suitable initial condition $\hat{w}_0$. In what follows, we require that the initial condition $\hat{w}_0$ is \emph{consistent}; see, for instance, \cite[p.\,456]{Hairer_Wanner}. This means that $\hat{w}_0$ satisfies the following conditions:
\begin{equation}\label{consistent_initial_conditions}
\left\{
\begin{array}{rl}
\vspace{0.25ex}
\mathcal{G}_u(u_0,t_0)\mathcal{F}(\hat{w}_0, t_0) = 0 \\
\vspace{0.25ex}
\mathcal{G}(u_0,t_0) = 0,
\end{array}
\right.
\end{equation}
where $\mathcal{G}_u$ is the derivative of $\mathcal{G}$ with respect to $u$. In this case, and only in this case, the index 2 DAE in \eqref{general_DAE} admits a locally unique solution; see, for example, the discussion in \cite[p.\,456]{Hairer_Wanner}.

\section{Runge--Kutta Discretization of the Navier--Stokes equations}\label{sec_3}
Given a domain $\Omega \subset \mathbb{R}^d$, $d=2,3$, we integrate the system of equations for $t\in[0,T]$:
\begin{equation}\label{Navier_Stokes_equation}
\left\{
\begin{array}{rl}
\vspace{0.25ex}
\frac{\partial u}{\partial t} - \nu \nabla^2 u + u \cdot \nabla u + \nabla p = f(\mathbf{x},t) & \quad \mathrm{in} \; \Omega \times (0,T), \\
\vspace{0.25ex}
\nabla \cdot u = 0 & \quad \mathrm{in} \; \Omega \times (0,T), \\
\vspace{0.25ex}
u(\mathbf{x},t) = g(\mathbf{x},t) & \quad \mathrm{on} \; \partial \Omega \times (0,T),\\
u(\mathbf{x},0) = u_0(\mathbf{x}) & \quad \mathrm{in} \; \Omega.

\end{array}
\right.
\end{equation}
Here, the \emph{force} function $f$, the boundary conditions $g$, and the initial condition $u_0$ are given. In what follows we assume  $u_0$ is solenoidal. 

We assume an inf-sup stable pair of finite elements for the spatial discretization of \eqref{Navier_Stokes_equation}. Denote the finite element  basis functions by $\{ \phi_i \}_{i=1}^{n_u}$ and $\{ \psi_i \}_{i=1}^{n_p}$, with $n_u,n_p \in \mathbb{N}$ the dimensions of the velocity and pressure spaces, respectively.

 To apply the framework of the previous section, we rewrite equations \eqref{Navier_Stokes_equation} as follows:
\begin{equation}\label{Navier_Stokes_equation_2}
\left\{
\begin{array}{rl}
\vspace{0.25ex}
\frac{\partial u}{\partial t} = \mathcal{F}(u, p, t) := f(\mathbf{x},t) + \nu \nabla^2 u - u \cdot \nabla u - \nabla {p} & \quad \mathrm{in} \; \Omega \times (0,T), \\
\vspace{0.25ex}
\mathcal{G}(u) := - \nabla \cdot u = 0 & \quad \mathrm{in} \; \Omega \times (0,T), \\
\vspace{0.25ex}
u(\mathbf{x},t) = g(\mathbf{x},t) & \quad \mathrm{on} \; \partial \Omega \times (0,T),\\
u(\mathbf{x},0) = u_0(\mathbf{x}) & \quad \mathrm{in} \; \Omega.
\end{array}
\right.
\end{equation}

For an index 2 DAE a Runge--Kutta time stepping technique is given by \eqref{DAE_RK_step}, where $\left\{Y^{u}_{i,n} \right\}_{i=1}^s$ and $\left\{Y^{p}_{i,n} \right\}_{i=1}^s$ solve the non-linear system in \eqref{stages_Runge_Kutta_2}, with $\mathcal{F}$ and $\mathcal{G}$ defined as in \eqref{Navier_Stokes_equation_2}. It is convenient to introduce the variables $w^{u}_{i,n}$ and $w^{p}_{i,n}$ as in \eqref{wu_wp}, where  $\left\{Y^{u}_{i,n} \right\}_{i=1}^s$ and $\left\{Y^{p}_{i,n} \right\}_{i=1}^s$ solve
\begin{equation}\label{RK_NS_non_linear_system}
\left\{
\begin{array}{ll}
\vspace{0.25ex}
Y^{u}_{i,n} = f(\mathbf{x}, t_n + c_i \Delta t) + \nu \nabla^2 w^{u}_{i,n} - w^{u}_{i,n} \cdot \nabla w^{u}_{i,n} - \nabla {w^{p}_{i,n}}, & \quad i=1,\ldots,s,\\
- \nabla \cdot w^{u}_{i,n} = 0, & \quad i=1,\ldots,s,
\end{array}
\right.
\end{equation}
with suitable boundary conditions. Since $\left\{Y^{u}_{i,n} \right\}_{i=1}^s$ represent an approximation of the time derivative of $u$, suitable boundary conditions on $\left\{Y^{u}_{i,n} \right\}_{i=1}^s$ are given by the time derivative of the boundary conditions on $u$.

We denote with $(\cdot, \cdot)$ the $L^2$-inner product in $\Omega$ and introduce the bilinear and trilinear forms
\begin{displaymath}
a_{\mathrm{NS}}(w, \phi) := \nu (\nabla w, \nabla \phi), \quad
b_{\mathrm{NS}}(w, \psi) := - (\nabla \cdot w, \psi), \quad
c_{\mathrm{NS}}(w; u, \phi) :=  (w \cdot \nabla u, \phi).
\end{displaymath}

Now we can write a finite element formulation of  \eqref{RK_NS_non_linear_system}:
\begin{equation}\label{RK_NS_weak_form}
\left\{
\begin{array}{ll}
\vspace{0.25ex}
(Y^{u}_{i,n}, \phi) + a_{\mathrm{NS}}(w^{u}_{i,n}, \phi) + c_{\mathrm{NS}}(w^{u}_{i,n}; \nabla w^{u}_{i,n}, \phi)\\
\phantom{(Y^{u}_{i,n}, \phi) +} + b_{\mathrm{NS}}(\phi, w^{p}_{i,n})= (f(\mathbf{x}, t_n + c_i \Delta t), \phi) , & \quad i=1,\ldots,s,\\
b_{\mathrm{NS}}(w^{u}_{i,n}, \psi) = 0, & \quad i=1,\ldots,s,
\end{array}
\right.
\end{equation}
where $\phi$ and $\psi$ are finite element basis functions.

The non-linearity of the trilinear forms $c_{\mathrm{NS}}(w^{u}_{i,n}; \nabla w^{u}_{i,n}, \phi)$, for $i=1,\ldots,s$, in \eqref{RK_NS_weak_form} requires us to solve for a linearization of \eqref{RK_NS_weak_form}. Our method of choice is a Newton linearization explained below.

Let $M_u$ and $M_p$ be the \emph{mass matrices} defined over the velocity and pressure spaces, respectively. Let $\mathbf{u}_n$ and $\mathbf{p}_n$ denote the coefficient vectors of the velocity and pressure at time $t_n = n \Delta t$, starting from $\mathbf{u}_0$ and $\mathbf{p}_0$, which are the coefficient vectors of the $L^2$-projections of the initial condition $(u_0, p_0)$ onto the finite element spaces. We can then write the Runge--Kutta method \eqref{DAE_RK_step} in matrix-vector form as follows:
\begin{equation}\label{RK_discretization_for_vp_Navier_Stokes_equations}
\begin{array}{ll}
\vspace{0.25ex}
M_u \mathbf{u}_{n+1} = M_u \mathbf{u}_n + \Delta t M_u \sum_{i=1}^{s} b_i \mathbf{Y}^{u}_{i,n}, & \quad n = 0, \ldots, n_t-1,\\
M_p \mathbf{p}_{n+1} = M_p \mathbf{p}_n + \Delta t M_p \sum_{i=1}^{s} b_i \mathbf{Y}^{p}_{i,n}, & \quad n = 0, \ldots, n_t-1.
\end{array}
\end{equation}
Here, the vectors $\mathbf{Y}^{u}_{i,n}$ and $\mathbf{Y}^{p}_{i,n}$, for $i=1,\ldots,s$, are suitable approximations of the stages $Y^{u}_{i,n}$ and $Y^{p}_{i,n}$ of the method between time $t_n$ and $t_{n+1}$, respectively, and solve a discretized version of the system in \eqref{RK_NS_weak_form}. In particular, if we denote with $\mathbf{w}^u_{i,n}$ (resp., $\mathbf{w}^{p}_{i,n}$) the vectors containing the discretization of $w^u_{i,n}$ (resp., of $w^{p}_{i,n}$), the discretized stages solve the following non-linear problem:
\begin{equation}\label{RK_stages_Navier_Stokes_equations}
\left\{
\begin{array}{ll}
\vspace{0.25ex}
M_u \mathbf{Y}^{u}_{i,n} + \sum_{l=1}^{n_u}\left( a_{\mathrm{NS}}(\mathbf{w}^u_{i,n}, \phi_l) + c_{\mathrm{NS}}(\mathbf{w}^u_{i,n}; \mathbf{w}^u_{i,n}, \phi_l) + b_\mathrm{NS}(\phi_l, \mathbf{w}^p_{i,n}) \right) = \mathbf{f}_{i,n}, \\
\sum_{l=1}^{n_p} b_\mathrm{NS}(\mathbf{w}^u_{i,n}, \psi_l) = 0,
\end{array}
\right.
\end{equation}
for $i=1,\ldots,s$ and $n=0,\ldots,n_t-1$.

In what follows, we denote with $K_u$ the \emph{vector-stiffness matrix}, $N_u(\mathbf{w}^u_{i,n})$ the \emph{vector-convection matrix} with the wind given by $\mathbf{w}^u_{i,n}$, where
\begin{displaymath}
N_u(\mathbf{w}^u_{i,n}) = \left\{n_{lj} \right\}, \quad n_{lj} = \int_\Omega (\mathbf{w}^u_{i,n} \cdot \nabla \phi_j) \phi_l \, \mathrm{d} \Omega,
\end{displaymath}
and $H_u(\mathbf{w}^u_{i,n})$ the matrix arising from second-order information on the convection term, namely,
\begin{displaymath}
H_u(\mathbf{w}^u_{i,n}) = \left\{h_{lj} \right\}, \quad h_{lj} = \int_\Omega (\phi_j \cdot \nabla \mathbf{w}^{u}_{i,n}) \phi_l \, \mathrm{d} \Omega.
\end{displaymath}
Further, we denote with $B$ the (negative) divergence matrix,
\begin{displaymath}
B = \left\{ - \int_\Omega \psi_l  \nabla \cdot \phi_j \, {\rm d}\Omega \right\}_{l = 1, \ldots, n_p}^{j = 1, \ldots, n_u}.
\end{displaymath}
Then, given an approximation $\mathbf{\hat{w}}^{(k)}_{i,n}=[\mathbf{w}^{u,(k)}_{i,n},\mathbf{w}^{p,(k)}_{i,n}]$ of $\mathbf{w}_{i,n} = [\mathbf{w}^{u}_{i,n},\mathbf{w}^{p}_{i,n}]$ and  recalling that $\mathbf{w}^{u}_{i,n}= \mathbf{u}_n + \Delta t \sum_{j=1}^{s} a_{i,j} \mathbf{Y}^{u}_{j,n}$ and $\mathbf{w}^{p}_{i,n}= \mathbf{p}_n + \Delta t \sum_{j=1}^{s} a_{i,j} \mathbf{Y}^{p}_{j,n}$, we linearize the equations in \eqref{RK_stages_Navier_Stokes_equations} as follows:
\begin{equation}\label{discretized_RK_stages_Navier_Stokes_equations}
\left\{
\begin{array}{ll}
\vspace{0.25ex}
M_u \mathbf{Y}^{u, (k+1)}_{i,n} + \Delta t \sum_{j=1}^{s} a_{i,j} L_u(\mathbf{w}^{u,(k)}_{i,n}) \mathbf{Y}^{u, (k+1)}_{j,n} + \Delta t \sum_{j=1}^{s} a_{i,j} B^\top \mathbf{Y}^{p, (k+1)}_{j,n} = \mathbf{\hat{f}}_{i,n}, \\
\Delta t \sum_{j=1}^{s} a_{i,j} B \mathbf{Y}^{u, (k+1)}_{j,n} = - B \mathbf{u}_n,
\end{array}
\right.
\end{equation}
for $i=1,\ldots,s$ and $n=0,\ldots,n_t-1$. In the equations above, we set $L_u(\mathbf{w}^{u,(k)}_{i,n}) := \nu K_u + N_u(\mathbf{w}^{u,(k)}_{i,n}) + H_u(\mathbf{w}^{u,(k)}_{i,n})$ and $\mathbf{\hat{f}}_{i,n} := \mathbf{f}_{i,n} - L_u(\mathbf{w}^{u,(k)}_{i,n}) \mathbf{u}_n - B^\top \mathbf{p}_n$.

The algebraic linear systems to be solved at each time step are given by
\begin{equation}\label{RK_NS_linear_system}
\underbrace{\left[
\begin{array}{cc}
\Phi^{(k)} & \Psi_1\\
\Psi_2 & 0
\end{array}
\right]}_{\mathcal{A}^{(k)}}
\left[
\begin{array}{l}
\mathbf{Y}^{u, (k+1)}_{1,n}\\
\vdots\\
\mathbf{Y}^{u, (k+1)}_{s,n}\\
\mathbf{Y}^{p, (k+1)}_{1,n}\\
\vdots\\
\mathbf{Y}^{p, (k+1)}_{s,n}
\end{array}
\right]
= 
\left[
\begin{array}{l}
\mathbf{b}^u_1\\
\vdots\\
\mathbf{b}^u_s\\
\mathbf{b}^p_1\\
\vdots\\
\mathbf{b}^p_s
\end{array}
\right],
\end{equation}
with the right-hand side accounting for the non-linear residual. Further, the $(1,1)$-block $\Phi^{(k)}$ of the matrix $\mathcal{A}^{(k)}$ is given by
\begin{displaymath}
\Phi^{(k)} =
\left[
\begin{array}{ccc}
M_u + \Delta t a_{1,1} L_u(\mathbf{w}^{u,(k)}_{1,n}) & \ldots & \Delta t a_{1,s} L_u(\mathbf{w}^{u,(k)}_{1,n})\\
\vdots & & \vdots\\
\Delta t a_{s,1} L_u(\mathbf{w}^{u,(k)}_{s,n}) & \ldots & M_u + \Delta t a_{s,s} L_u(\mathbf{w}^{u,(k)}_{s,n})
\end{array}
\right],
\end{displaymath}
while $\Psi_1 = \Delta t A_{\mathrm{RK}}\otimes B^\top$ and $\Psi_2 = \Delta t A_{\mathrm{RK}}\otimes B$.

We conclude this section by discussing how to find consistent initial conditions. Since $u_0$ is solenoidal, the second condition in \eqref{consistent_initial_conditions} is satisfied. Furthermore, from the definition of $\mathcal{G}(\cdot)$ in \eqref{Navier_Stokes_equation_2}, we have $\mathcal{G}_{u}(\cdot) = (-\nabla \cdot \,)$. Therefore, for the incompressible Navier--Stokes equations, the first condition in \eqref{consistent_initial_conditions} takes the following form:
\begin{align*}
\mathcal{G}_{u}(u_0,t_0)\mathcal{F}(u_0, p_0, t_0) & = - \nabla \cdot( f(\mathbf{x},t_0) + \nu \nabla^2 u_0 - u_0 \cdot \nabla u_0 - \nabla {p}_0) \\
& = \nabla \cdot \nabla p_0 - \nabla \cdot (f(\mathbf{x},t_0) + \nu \nabla^2 u_0 - u_0 \cdot \nabla u_0)\\
& = \nabla^2 p_0 - \nabla \cdot (f(\mathbf{x},t_0) + \nu \nabla^2 u_0 - u_0 \cdot \nabla u_0) = 0,
\end{align*}
which implies that
\begin{displaymath}
- \nabla^2 p_0 = - \nabla \cdot (f(\mathbf{x},t_0) + \nu \nabla^2 u_0 - u_0 \cdot \nabla u_0).
\end{displaymath}
In practice, a consistent initial condition can be found by employing a Laplace solver on the pressure space.

\section{Preconditioner}\label{sec_4}
This work aims to develop a fast and robust preconditioner for the sequence of Newton linearizations of \eqref{RK_stages_Navier_Stokes_equations} arising at each time step. Each system has the form \eqref{RK_NS_linear_system}, a \emph{generalized saddle-point problem} as defined in \cite{Benzi_Golub_Liesen}. We begin with a brief overview of preconditioning such systems, then introduce the strategy used in our setting. For simplicity, we omit the superscript indicating the nonlinear iteration in what follows.

Suppose we aim to solve a generalized saddle-point problem of the form
\begin{equation}\label{generalized_saddle_point_problem}
\underbrace{\left[
\begin{array}{cc}
\Phi & \Psi_1\\
\Psi_2 & -\Theta
\end{array}
\right]}_{\mathcal{A}}
\left[
\begin{array}{l}
\mathbf{Y}_1\\
\mathbf{Y}_2
\end{array}
\right]
=
\left[
\begin{array}{l}
\mathbf{b}_1\\
\mathbf{b}_2
\end{array}
\right],
\end{equation}
with $\mathcal{A}$ invertible. If the $(1,1)$ block $\Phi$ and the Schur complement $S := \Theta + \Psi_2 \Phi^{-1} \Psi_1$ are invertible, then, as shown in \cite{Murphy_Golub_Wathen}, an ideal preconditioner for $\mathcal{A}$ is given by
\begin{equation}\label{optimal_preconditioner}
\mathcal{P} = \left[
\begin{array}{cc}
\Phi & \Psi_1\\
0 & -S
\end{array}
\right].
\end{equation}
In this case, the preconditioned matrix has all eigenvalues equal to 1, and its minimal polynomial has degree two. Consequently, a Krylov subspace method would converge in at most two iterations in exact arithmetic\footnote{A similar argument is valid for the block-lower triangular version of $\mathcal{P}$}. Since the preconditioner $\mathcal{P}$ is non-symmetric, GMRES is the appropriate solver \cite{Saad_Schultz}.

In practical applications, rather than solving for the matrix $\mathcal{P}$ above, one employs a suitable approximation $\widehat{\mathcal{P}}$ of it of the form
\begin{equation}\label{approx_preconditioner}
\widehat{\mathcal{P}}=\left[
\begin{array}{cc}
\widehat{\Phi} & \Psi_1\\
0 & -\widehat{S}
\end{array}
\right].
\end{equation}

We now move to describing of the preconditioning strategy employed in this work. Let us consider the generalized saddle-point problem in \eqref{RK_NS_linear_system}. Given a suitable invertible matrix $\mathcal{W}$ and $\gamma>0$, if we set $
\mathbf{Y}^u_n=[\mathbf{Y}^u_{1,n}, \ldots, \mathbf{Y}^u_{s,n}]^\top$ and $
\mathbf{Y}^p_n=[\mathbf{Y}^p_{1,n}, \ldots, \mathbf{Y}^p_{s,n}]^\top$, rather than considering the solution of \eqref{RK_NS_linear_system}, we consider an equivalent system of the form
\begin{equation}\label{augmented_system}
\left[
\begin{array}{cc}
\Phi + \gamma \Psi_1 \mathcal{W}^{-1} \Psi_2 & \Psi_1\\
\Psi_2 & 0
\end{array}
\right]
\left[
\begin{array}{l}
\mathbf{Y}^u\\
\mathbf{Y}^p
\end{array}
\right]
= 
\left[
\begin{array}{l}
\mathbf{\hat{b}}^u\\
\mathbf{b}^p
\end{array}
\right],
\end{equation}
where for the right-hand side we set $\mathbf{\hat{b}}^u = \mathbf{b}^u + \gamma \Psi_1 \mathcal{W}^{-1} \mathbf{b}^p$, with $\mathbf{b}^u=[\mathbf{b}^u_1, \ldots, \mathbf{b}^u_s]^\top$ and $\mathbf{b}^p = [\mathbf{b}^p_1, \ldots, \mathbf{b}^p_s]^\top$. We employ as matrix $\mathcal{W}$ the following matrix:
\begin{equation}\label{mathcal_W}
\mathcal{W} = \Delta t A_{\mathrm{RK}} \otimes W,
\end{equation}
for a suitable $W$. In the context of incompressible Navier--Stokes equations, the matrix $W$ is chosen as either the pressure mass matrix $M_p$ or its diagonal. We will specify this in our numerical tests below.

Then, the preconditioner $\mathcal{P}$ defined in \eqref{optimal_preconditioner} specifies as follows:
\begin{equation}\label{AL_preconditioner}
\mathcal{P}=\left[
\begin{array}{cc}
\Phi + \gamma \Psi_1 \mathcal{W}^{-1} \Psi_2 & \Psi_1\\
0 & -S_\gamma
\end{array}
\right],
\end{equation}
where $S_\gamma = \Psi_2 ( \Phi + \gamma \Psi_1 \mathcal{W}^{-1} \Psi_2)^{-1} \Psi_1$.

We now proceed with finding suitable cheap approximations of the diagonal blocks of $\mathcal{P}$, starting from the $(1,1)$-block $\Phi + \gamma \Psi_1 \mathcal{W}^{-1} \Psi_2$ and then moving to the approximation of the Schur complement $S_\gamma$.

The $(1,1)$-block $\Phi + \gamma \Psi_1 \mathcal{W}^{-1} \Psi_2$ of the preconditioner $\mathcal{P}$ has a complex structure. It is an $s \times s$ block matrix, where each block-row contains a spatial discretization of a differential operator that can be viewed as a ``perturbation" of a convection--diffusion operator with a time-dependent wind field.
In recent years, several solvers have been proposed for the discretization of time-dependent convection--diffusion problems using Runge--Kutta methods in time; see, for example, \cite{Axelsson_Dravins_Neytcheva, Kirby, Mardal_Nilssen_Staff, Rana_Howle_Long_Meek_Milestone, Southworth_Krzysik_Pazner_Sterck_2, VanLent_Vanderwalle}. These works exploit the Kronecker product structure of the resulting systems to design and analyze efficient solvers. Unfortunately, the discretized system in our setting lacks this Kronecker structure, which limits the applicability of those methods to heuristic use only.
To the best of our knowledge, the only work addressing problems of this form is \cite{AbuLabdeh_MachLachlan_Farrell}, where the authors developed a monolithic multigrid solver for the nonstationary incompressible Stokes equations discretized in time with a Runge--Kutta method. They extended this approach to the incompressible Navier--Stokes equations and further to magnetohydrodynamic problems.
Although a monolithic treatment of the $(1,1)$-block is possible in principle, we instead adopt a block Gauss--Seidel approximation, as proposed in \cite{Staff_Mardal_Nilssen}. A detailed comparison of the aforementioned techniques for approximating the $(1,1)$-block will be the subject of future work.

We now focus on the Schur complement $S_\gamma$. To approximate the Schur complement $S_\gamma$, we employ the Sherman--Morrison--Woodbury formula. Assuming that $\Psi_1 \mathcal{W}^{-1} \Psi_2$ is invertible, we can write:
\begin{displaymath}
S_\gamma^{-1} = \gamma \mathcal{W}^{-1} + (\Psi_2\Phi^{-1}\Psi_1)^{-1}.
\end{displaymath}
Assuming for the moment that this inverse exists, we now discuss how to approximate each of its terms.

From \eqref{mathcal_W}, using well-known properties of the Kronecker product and recalling that $A_{\mathrm{RK}}$ is invertible, we obtain:
\begin{displaymath}
\mathcal{W}^{-1} = \dfrac{1}{\Delta t} A^{-1}_{\mathrm{RK}} \otimes W^{-1}.
\end{displaymath}
In our numerical tests, we use a matrix $W$ that is block diagonal, with each block being a mass matrix on the pressure space or its diagonal, so it is easily invertible.

Next, we consider the matrix $\Psi_2\Phi^{-1}\Psi_1$. Again using properties of the Kronecker product, we can write:
\begin{align*}
\Psi_2\Phi^{-1}\Psi_1 & = \Delta t^2 (A_{\mathrm{RK}} \otimes B) \Phi^{-1}(A_{\mathrm{RK}} \otimes B^\top)\\
 & = \Delta t^2 (A_{\mathrm{RK}} \otimes I_p)\, \underbrace{(I_s \otimes B) \Phi^{-1}(I_s \otimes B^\top)}_{S_{\mathrm{int}}} \,(A_{\mathrm{RK}} \otimes I_p),
\end{align*}
where $I_p$ is the identity matrix on the pressure space. From here, we have
\begin{displaymath}
(\Psi_2\Phi^{-1}\Psi_1)^{-1} = \Delta t^{-2} (A^{-1}_{\mathrm{RK}} \otimes I_p) S^{-1}_{\mathrm{int}} (A^{-1}_{\mathrm{RK}} \otimes I_p).
\end{displaymath}
The expression above shows that a low-cost approximation of $S_\gamma$ can be obtained if a suitable approximation of $S^{-1}_{\mathrm{int}}$ is available. We adopt the following heuristic approximation:
\begin{equation}\label{S_tilde}
I_s \otimes K_p^{-1} + \nu \Delta t A_{\mathrm{RK}} \otimes M_p^{-1} =:\widetilde{S}^{-1}_{\mathrm{int}} \approx S^{-1}_{\mathrm{int}}.
\end{equation}
This approximation was proposed in \cite{Leveque_Bergamaschi_Martinez_Pearson} for the Stokes equations, in the limit of very viscous flow.

Finally, we define the approximation $\widetilde{S}^{-1}_\gamma$ of $S^{-1}_\gamma$ used in our numerical tests:
\begin{equation}\label{Schur_approx}
\widetilde{S}^{-1}_\gamma = \gamma \mathcal{W}^{-1} + \Delta t^{-2} (A^{-1}_{\mathrm{RK}} \otimes I_p) \widetilde{S}^{-1}_{\mathrm{int}} (A^{-1}_{\mathrm{RK}} \otimes I_p).
\end{equation}

We note that the augmented Lagrangian preconditioner described above was derived from an algebraic standpoint. However, in the case of the incompressible Navier--Stokes equations, a similar structure for the linearized problem in \eqref{augmented_system} can be obtained by adding a grad-div stabilization term~\cite{gradDiv} to the differential operator. {In fact, in this case the linearized system would have the form as in \eqref{augmented_system}, with the block $\gamma \Psi_1 \mathcal{W}^{-1} \Psi_2$ replaced by $\gamma A_{\mathrm{RK}} \otimes R$, where $R$ is the discrete grad-div operator}. For this reason, the preconditioner introduced in this section is also applicable to grad-div stabilized discretizations of the incompressible Navier--Stokes equations, as we will show in the numerical results.

\section{Numerical results}\label{sec_5}
We now assess the efficiency and robustness of the augmented Lagrangian preconditioner derived in this work. In all our tests $d=2$. All CPU times below are reported in seconds.

\subsection{Experiments with exact block-diagonal solvers}
We first present numerical results using the augmented Lagrangian approach from an algebraic perspective, employing only Radau IIA methods. This choice is motivated by their suitability for index-2 DAEs: Radau IIA methods are L-stable and do not suffer from order reduction, unlike Gauss or Lobatto IIIC methods (see, e.g., \cite[Table 4.1, p. 504]{Hairer_Wanner}).

We begin by studying how the augmented Lagrangian preconditioner affects the discretization error, then examine its dependence on problem parameters when Radau IIA methods are used in time. The influence of the Runge--Kutta method itself is considered in the next section. All tests were run in MATLAB R2018b on an Ubuntu 18.04.1 LTS system with a 1.70 GHz Intel quad-core i5 processor and 8 GB RAM.

We use inf--sup stable Taylor--Hood $\mathbf{Q}_2$--$\mathbf{Q}_1$ elements for spatial discretization of the Navier--Stokes equations. For a refinement level $l$, we use a uniform mesh of size $h = 2^{-1 - l}$ for the velocity and $h = 2^{-l}$ for the pressure in each spatial direction. Letting $q_{\mathrm{FE}}$ be the finite element order and $q_{\mathrm{RK}}$ the Runge--Kutta order, we set the number of time steps $n_t$ so that $\Delta t = T / n_t \leq h^{q_{\mathrm{FE}} / q_{\mathrm{RK}}}$.

At high Reynolds numbers, the problem becomes convection-dominated and requires stabilization. We adopt the Local Projection Stabilization (LPS) method as described in \cite{Becker_Braack, Becker_Vexler, Braack_Burman}. For convergence analysis of LPS applied to the Oseen problem, see \cite{Matthies_Tobiska}. Alternative stabilization techniques can be found in \cite{Brooks_Hughes, Franca_Frey, Johnson_Saranen, Tobiska_Lube}. Throughout, $\text{Id}$ denotes the identity operator, and the union of elements sharing a vertex is referred to as a patch.

Given a convection field $v$, the LPS stabilization is defined using the fluctuation operator $\kappa_h = \text{Id} - \pi_h$, where $\pi_h$ is the $L^2$-orthogonal (discontinuous) projection onto a patch-wise finite-dimensional space $\overline{V}$. Assuming that the triangulation of $\Omega$ consists of $N_{\rm patch}$ patches, the stabilization matrix $Q_u$ is given by
\begin{displaymath}
Q_u = \left\{ \sum_{m=1}^{N_{\rm patch}} \int_{\mathtt{P}_m} \delta_m \kappa_h(v \cdot \nabla \phi_i) \cdot \kappa_h(v \cdot \nabla \phi_j) \, {\rm d}\mathtt{P}_m \right\}_{i,j=1}^{n_v},
\end{displaymath}
where $\delta_m \ge 0$ is a stabilization parameter {and $\mathtt{P}_m$ is an elements patch}. Following \cite{Becker_Braack}, we take
\begin{displaymath}
\left.\pi_h(q)\right|_{\mathtt{P}_m} = \dfrac{1}{|\mathtt{P}_m|} \int_{\mathtt{P}_m} q \: \mathrm{d} \mathtt{P}_m, \quad \forall q \in L^2(\Omega),
\end{displaymath}
where $|\mathtt{P}_m|$ is the (Lebesgue) measure of an elements patch $\mathtt{P}_m$. 
Further, as in \cite[p.\,253]{Elman_Silvester_Wathen} we set
\begin{displaymath}
\delta_m =
\left\{
\begin{array}{lcl}
\vspace{1ex}
\dfrac{h_m}{2 \| v(\mathbf{x}_m) \| } \left(1 - \dfrac{1}{Pe_m}\right) & \quad \mathrm{if} \: Pe_m > 1, \\
0 & \quad \mathrm{if} \: Pe_m \leq 1.
\end{array}
\right.
\end{displaymath}
where $\| v(\mathbf{x}_m) \|$ is the $\ell_2$-norm of $v$ at the patch centroid $\mathbf{x}_m$, $h_m$ is a characteristic patch length in the flow direction, and $Pe_m = \| v(\mathbf{x}_m) \| h_m / (2 \nu)$ is the patch Péclet number.

We take the diagonal of the pressure mass matrix as $W$ in \eqref{mathcal_W}. Applying the inverse of $\widetilde{S}_\mathrm{int}$ (from \eqref{S_tilde}) requires repeated inversion of the pressure mass and stiffness matrices. The mass matrix is approximately inverted using 20 steps of Chebyshev semi-iteration \cite{GolubVargaI, GolubVargaII, Wathen_Rees} with Jacobi preconditioning. The pressure stiffness matrix has a kernel due to the hydrostatic pressure mode; to handle this, we fix the pressure at one node and approximate its inverse using two V-cycles of the \texttt{HSL\_MI20} solver \cite{HSL_MI20}.

For the inversion of the diagonal blocks in the Gauss--Seidel method applied to the $(1,1)$-block, we use MATLAB's backslash operator. This is our choice here in lieu of a tailored geometric multigrid solver for the augmented system. A nonlinear multigrid approach is discussed in the next section as an alternative. Consequently, we use the flexible GMRES algorithm \cite{Saad} as the outer solver in all experiments, with a stopping tolerance of $10^{-6}$ on the residual. Newton’s method is run until the relative nonlinear residual is reduced by $10^{-5}$.

\subsubsection{Accuracy test}
In this section, we assess the discretization errors of the recovered finite element solutions when using an $s$-stage Radau IIA method for time integration. Since the augmentation alters the system residual, we investigate whether it affects the finite element accuracy.

We consider the domain $\Omega = (0,1)^2$, final time $T = 2$, and solve problem \eqref{Navier_Stokes_equation} with $\nu = \frac{1}{50}$, using the exact solution:
\begin{displaymath}
\begin{array}{l}
\vspace{0.2ex}
u(\mathbf{x},t) = \frac{1}{2}e^{T - t} [\sin(\pi x_1) \cos(\pi x_2), - \cos(\pi x_1) \sin(\pi x_2)]^\top,\\
p(\mathbf{x},t) = \text{constant},
\end{array}
\end{displaymath}
with the forcing term $f$, as well as boundary and initial conditions, derived accordingly.
We evaluate the velocity error in the $L^\infty(\mathcal{H}_0^1(\Omega)^d)$ norm for the velocity and in the $L^\infty(L^2(\Omega))$ norm for the pressure, defined respectively as follows:
\begin{displaymath}
\begin{array}{c}
\displaystyle u_{\text{error}} = \max_{n} \left[ (\mathbf{u}_n - \mathbf{u}_n^{\mathrm{sol}})^\top K_u (\mathbf{u}_n - \mathbf{u}_n^{\mathrm{sol}}) \right]^{1/2}, \\
\displaystyle p_{\text{error}} = \max_{n} \left[ (\mathbf{p}_n - \mathbf{p}_n^{\mathrm{sol}})^\top M_p (\mathbf{p}_n - \mathbf{p}_n^{\mathrm{sol}}) \right]^{1/2},
\end{array}
\end{displaymath}
where $\mathbf{u}_n^{\mathrm{sol}}$ and $\mathbf{p}_n^{\mathrm{sol}}$ denote the discretized exact solutions for $u$ and $p$ at time $t_n$, respectively.

Table \ref{table_convergence_Navier_Stokes_equations} reports the average number of FGMRES iterations ($\texttt{it}$), the average CPU time per linear iteration, and the discretization errors for velocity and pressure. Table \ref{table_number_newton_iter} presents the average number of Newton iterations per time step, along with the problem size (degrees of freedom, DoF).

\begin{table}[!ht]
\caption{Average FGMRES iterations, average CPU times, and resulting errors in $u$ and $p$, for Radau IIA methods with $\gamma=1$ and $\gamma=100$.}\label{table_convergence_Navier_Stokes_equations}
\begin{center}
\begin{footnotesize}
{\begin{tabular}{|c|c||c|c|c|c||c|c|c|c|}
\hline
\multicolumn{1}{|c}{} & \multicolumn{1}{|c||}{} & \multicolumn{4}{c||}{$\gamma=1$} & \multicolumn{4}{c|}{$\gamma=100$} \\
\cline{3-10}
$s$ & $l$ & $\texttt{it}$ & CPU & $u_{\text{error}}$ & $p_{\text{error}}$ & $\texttt{it}$ & CPU & $u_{\text{error}}$ & $p_{\text{error}}$ \\
\hline
\hline
\multirow{5}{*}{2} & $2$ & 13 & 0.02 & 1.46e+00 & 2.74e-01 & 8 & 0.01 & 1.46e+00 & 2.74e-01 \\
 & $3$ & 13 & 0.16 & 7.50e-01 & 9.30e-02 & 8 & 0.10 & 7.50e-01 & 9.30e-02 \\
 & $4$ & 12 & 1.07 & 7.24e-02 & 3.01e-03 & 7 & 0.58 & 7.24e-02 & 3.01e-03 \\
 & $5$ & 11 & 7.26 & 1.13e-03 & 9.60e-06 & 7 & 4.10 & 1.13e-03 & 7.95e-06 \\
 & $6$ & 10 & 36.9 & 1.50e-04 & 3.42e-06 & 6 & 22.3 & 1.50e-04 & 1.75e-06 \\
\hline
\hline
\multirow{5}{*}{4} & $2$ & 18 & 0.05 & 1.03e+00 & 2.16e-01 & 13 & 0.03 & 1.03e+00 & 2.16e-01 \\
 & $3$ & 21 & 0.59 & 3.80e-01 & 9.46e-02 & 13 & 0.32 & 3.80e-01 & 9.46e-02 \\
 & $4$ & 21 & 4.05 & 9.08e-03 & 7.54e-04 & 12 & 2.10 & 9.08e-03 & 7.53e-04 \\
 & $5$ & 20 & 25.1 & 6.99e-04 & 8.25e-06 & 10 & 12.9 & 6.99e-04 & 2.08e-06 \\
 & $6$ & 19 & 136 & 9.53e-05 & 7.80e-07 & 10 & 70.0 & 9.51e-05 & 6.11e-07 \\
\hline
\end{tabular}}
\end{footnotesize}
\end{center}
\end{table}

\begin{table}[!ht]
\caption{Degrees of freedom (DoF) and average number of Newton iterations per time step, with $\gamma=1$ and $\gamma=100$.}\label{table_number_newton_iter}
\begin{small}
\begin{center}
\renewcommand{\arraystretch}{1.2}
\begin{tabular}{|c||c|c|c||c|c|c|}
\hline
 & \multicolumn{3}{c||}{$s=2$} & \multicolumn{3}{c|}{$s=4$} \\
\cline{2-7}
$l$ & DoF & $\gamma=1$ & $\gamma=100$ & DoF & $\gamma=1$ & $\gamma=100$ \\
\hline
\hline
$2$ & 246 & 4 & 4 & 492 & 4 & 5 \\
\hline
$3$ & 1062 & 3 & 4 & 2124 & 4 & 5 \\
\hline
$4$ & 4422 & 2 & 5 & 8844 & 4 & 6 \\
\hline
$5$ & 18,054 & 2 & 5 & 36,108 & 4 & 6 \\
\hline
$6$ & 72,966 & 3 & 6 & 145,932 & 4 & 8 \\
\hline
\end{tabular}
\end{center}
\end{small}
\end{table}

Table \ref{table_convergence_Navier_Stokes_equations} shows that the average number of FGMRES iterations remains nearly constant under mesh refinement. CPU times scale super-linearly with system size, which is expected due to the exact inversion of the augmented diagonal blocks in the block Gauss--Seidel method. Increasing $\gamma$ reduces linear iterations, while velocity errors remain unaffected and pressure errors show mild sensitivity to $\gamma$ on finer meshes. From Table \ref{table_number_newton_iter} we observe that the average number of Newton iterations per time step is stable, with a slight increase for larger $\gamma$.

\subsubsection{Lid-driven cavity}\label{sec_5_1_2}
This section studies the dependence of the proposed preconditioner on the parameters of both the problem and the numerical method. We continue using $s$-stage Radau IIA methods for time discretization. Specifically, we investigate how the performance of the preconditioner depends on the problem parameters by varying the mesh size $h$, the viscosity $\nu$, the parameter $\gamma$, and the number of stages $s$. Recall that we choose the time step $\Delta t$ such that $\Delta t \leq h^{q_{\mathrm{FE}}/q_{\mathrm{RK}}}$. In the next section, we analyze how the preconditioner depends on the class of Runge--Kutta methods used.

For this test, we set $\Omega = (-1,1)^2$ and $T = 2$. Starting from the initial condition $u(\mathbf{x},0) = [0,0]^\top$, the flow is governed by
\eqref{Navier_Stokes_equation} with $f(\mathbf{x}, t)=[0,0]^\top$ and the following boundary conditions:
\begin{displaymath}
g(\mathbf{x},t)=
\left\{
\begin{array}{cl}
\left[t,0\right]^\top & \mathrm{on} \: \partial \Omega_1 \times (0,1),\\
\left[1,0\right]^\top & \mathrm{on} \: \partial \Omega_1 \times [1,T),\\
\left[0,0\right]^\top & \mathrm{on} \: (\partial \Omega \setminus \partial\Omega_1) \times (0,T),
\end{array}
\right.
\end{displaymath}
with $\partial \Omega_1 := \left(-1,1 \right)\times \left\{1\right\}$.

In Tables \ref{Navier_Stokes_equations_table_Radau_Matlab_gamma1}--\ref{Navier_Stokes_equations_table_Radau_Matlab_gamma100}, we report the average number of FGMRES iterations (\texttt{it}), the average number of Newton iterations per time step ($\mathrm{Nit}$), and the average CPU time per linear solve. Additionally, Table \ref{Navier_Stokes_equations_table_Radau_Matlab_gamma1} includes the dimension (DoF) of the corresponding system to be solved.

\begin{table}[!ht]
\caption{Average FGMRES iterations ($\texttt{it}$), average Newton iterations ($\mathrm{Nit}$) per time step, average CPU times per linear iteration, and degrees of freedom for the block augmented Lagrangian preconditioner with $\gamma=1$, for $\nu=\frac{1}{100}$, $\frac{1}{250}$, and $\frac{1}{500}$ and a range of $l$, for Radau IIA methods.}\label{Navier_Stokes_equations_table_Radau_Matlab_gamma1}
\begin{footnotesize}
\begin{center}
\renewcommand{\arraystretch}{1.2}
{\begin{tabular}{|c||c||c||c|c|c||c|c|c||c|c|c|}
\hline
\multicolumn{1}{|c||}{} &\multicolumn{1}{c||}{} & \multicolumn{1}{c||}{} & \multicolumn{3}{c||}{$\nu=\frac{1}{100}$} & \multicolumn{3}{c||}{$\nu=\frac{1}{250}$} & \multicolumn{3}{c|}{$\nu=\frac{1}{500}$}\\
\cline{4-12}
$s$ & $l$ & DoF & $\texttt{it}$ & $\mathrm{Nit}$ & CPU & $\texttt{it}$ & $\mathrm{Nit}$ & CPU & $\texttt{it}$ & $\mathrm{Nit}$ & CPU \\
\hline
\hline
\multirow{4}{*}{2}
& $3$ & 1062 & 13 & 4 & 0.17 & 13 & 6 & 0.16 & 13 & 8 & 0.16 \\
& $4$ & 4422 & 13 & 3 & 1.08 & 12 & 4 & 1.04 & 13 & 7 & 1.07 \\
& $5$ & 18,054 & 12 & 3 & 6.68 & 12 & 3 & 6.61 & 12 & 4 & 6.75 \\
& $6$ & 72,966 & 12 & 3 & 41 & 12 & 3 & 39.3 & 12 & 4 & 41.3 \\
\hline
\multirow{4}{*}{3}
& $3$ & 1593 & 16 & 4 & 0.3 & 17 & 6 & 0.32 & 16 & 8 & 0.31 \\
& $4$ & 6633 & 16 & 3 & 2.17 & 17 & 5 & 2.17 & 17 & 8 & 2.24 \\
& $5$ & 27,081 & 16 & 3 & 14 & 16 & 4 & 13.5 & 16 & 4 & 14.2 \\
& $6$ & 109,449 & 16 & 3 & 85.5 & 16 & 4 & 84.3 & 16 & 4 & 90.9 \\
\hline
\multirow{4}{*}{4}
& $3$ & 2124 & 18 & 5 & 0.48 & 18 & 6 & 0.47 & 20 & 8 & 0.51 \\
& $4$ & 8844 & 20 & 4 & 3.51 & 20 & 5 & 3.54 & 20 & 7 & 3.55 \\
& $5$ & 36,108 & 21 & 4 & 24.5 & 20 & 4 & 23 & 21 & 5 & 25.2 \\
& $6$ & 145,932 & 20 & 3 & 145 & 20 & 4 & 144 & 20 & 4 & 153 \\
\hline
\multirow{4}{*}{5}
& $3$ & 2655 & 22 & 5 & 0.74 & 24 & 7 & 0.8 & 25 & 8 & 0.86 \\
& $4$ & 11,055 & 24 & 4 & 5.88 & 24 & 5 & 5.61 & 25 & 8 & 5.69\\
& $5$ & 45,135 & 25 & 4 & 36.6 & 24 & 4 & 35.9 & 27 & 5 & 40.8 \\
& $6$ & 182,415 & 25 & 4 & 232 & 25 & 4 & 235 & 25 & 4 & 247 \\
\hline
\end{tabular}}
\end{center}
\end{footnotesize}
\end{table}

\begin{table}[!ht]
\caption{Average FGMRES iterations ($\texttt{it}$), average Newton iterations ($\mathrm{Nit}$) per time step, and average CPU times per linear iteration for the block augmented Lagrangian preconditioner with $\gamma=10$, for $\nu=\frac{1}{100}$, $\frac{1}{250}$, and $\frac{1}{500}$ and a range of $l$, for Radau IIA methods.}\label{Navier_Stokes_equations_table_Radau_Matlab_gamma10}
\begin{footnotesize}
\begin{center}
\renewcommand{\arraystretch}{1.2}
{\begin{tabular}{|c||c||c|c|c||c|c|c||c|c|c|}
\hline
\multicolumn{1}{|c||}{} &\multicolumn{1}{c||}{} & \multicolumn{3}{c||}{$\nu=\frac{1}{100}$} & \multicolumn{3}{c||}{$\nu=\frac{1}{250}$} & \multicolumn{3}{c|}{$\nu=\frac{1}{500}$}\\
\cline{3-11}
$s$ & $l$ & $\texttt{it}$ & $\mathrm{Nit}$ & CPU & $\texttt{it}$ & $\mathrm{Nit}$ & CPU & $\texttt{it}$ & $\mathrm{Nit}$ & CPU \\
\hline
\hline
\multirow{4}{*}{2}
& $3$ & 10 & 4 & 0.13 & 10 & 6 & 0.13 & 10 & 8 & 0.13 \\
& $4$ & 11 & 4 & 0.83 & 10 & 4 & 0.82 & 10 & 7 & 0.83 \\
& $5$ & 10 & 4 & 5.58 & 10 & 4 & 5.5 & 10 & 4 & 5.47 \\
& $6$ & 10 & 4 & 35.1 & 10 & 4 & 34.4 & 10 & 4 & 36.4 \\
\hline
\multirow{4}{*}{3}
& $3$ & 15 & 5 & 0.28 & 16 & 7 & 0.29 & 16 & 8 & 0.29 \\
& $4$ & 15 & 4 & 1.88 & 15 & 5 & 1.97 & 16 & 8 & 2.07 \\
& $5$ & 14 & 4 & 11.4 & 14 & 4 & 11.4 & 15 & 4 & 12.6 \\
& $6$ & 13 & 4 & 71.7 & 14 & 4 & 72 & 14 & 4 & 77.9 \\
\hline
\multirow{4}{*}{4}
& $3$ & 19 & 5 & 0.47 & 19 & 6 & 0.47 & 20 & 8 & 0.51 \\
& $4$ & 18 & 4 & 3.14 & 19 & 5 & 3.33 & 19 & 7 & 3.34\\
& $5$ & 18 & 4 & 20.4 & 18 & 4 & 20.4 & 20 & 5 & 23.5 \\
& $6$ & 18 & 4 & 129 & 18 & 4 & 127 & 18 & 4 & 139 \\
\hline
\multirow{4}{*}{5}
& $3$ & 24 & 5 & 0.75 & 25 & 7 & 0.84 & 26 & 8 & 0.89 \\
& $4$ & 24 & 4 & 5.49 & 24 & 5 & 5.47 & 24 & 8 & 5.5 \\
& $5$ & 24 & 4 & 34.3 & 23 & 4 & 33.6 & 25 & 5 & 38.1 \\
& $6$ & 23 & 4 & 211 & 22 & 4 & 208 & 23 & 4 & 222 \\
\hline
\end{tabular}}
\end{center}
\end{footnotesize}
\end{table}

\begin{table}[!ht]
\caption{Average FGMRES iterations ($\texttt{it}$), average Newton iterations ($\mathrm{Nit}$) per time step, and average CPU times per linear iteration for the block augmented Lagrangian preconditioner with $\gamma=100$, for $\nu=\frac{1}{100}$, $\frac{1}{250}$, and $\frac{1}{500}$ and a range of $l$, for Radau IIA methods.}\label{Navier_Stokes_equations_table_Radau_Matlab_gamma100}
\begin{footnotesize}
\begin{center}
\renewcommand{\arraystretch}{1.2}
{\begin{tabular}{|c||c||c|c|c||c|c|c||c|c|c|}
\hline
\multicolumn{1}{|c||}{} & \multicolumn{1}{c||}{} & \multicolumn{3}{c||}{$\nu=\frac{1}{100}$} & \multicolumn{3}{c||}{$\nu=\frac{1}{250}$} & \multicolumn{3}{c|}{$\nu=\frac{1}{500}$}\\
\cline{3-11}
$s$ & $l$ & $\texttt{it}$ & $\mathrm{Nit}$ & CPU & $\texttt{it}$ & $\mathrm{Nit}$ & CPU & $\texttt{it}$ & $\mathrm{Nit}$ & CPU \\
\hline
\hline
\multirow{4}{*}{2}
& $3$ & 8 & 5 & 0.11 & 8 & 6 & 0.1 & 8 & 8 & 0.1 \\
& $4$ & 8 & 4 & 0.61 & 8 & 5 & 0.64 & 8 & 7 & 0.67 \\
& $5$ & 8 & 4 & 4.09 & 8 & 5 & 4.04 & 8 & 5 & 4.3 \\
& $6$ & 7 & 5 & 26.2 & 7 & 5 & 26.2 & 7 & 5 & 27.6 \\
\hline
\multirow{4}{*}{3}
& $3$ & 12 & 6 & 0.21 & 12 & 8 & 0.23 & 12 & 8 & 0.22 \\
& $4$ & 10 & 5 & 1.24 & 11 & 6 & 1.43 & 13 & 8 & 1.65 \\
& $5$ & 10 & 5 & 8.2 & 10 & 5 & 7.99 & 11 & 6 & 8.87 \\
& $6$ & 10 & 5 & 51.7 & 9 & 6 & 47.9 & 9 & 6 & 51.5 \\
\hline
\multirow{4}{*}{4}
& $3$ & 15 & 6 & 0.39 & 15 & 7 & 0.38 & 16 & 8 & 0.5 \\
& $4$ & 14 & 5 & 2.6 & 15 & 6 & 2.68 & 16 & 8 & 2.75 \\
& $5$ & 14 & 5 & 16.4 & 14 & 5 & 15.7 & 15 & 6 & 17.6 \\
& $6$ & 14 & 5 & 103 & 14 & 5 & 99 & 14 & 5 & 102 \\
\hline
\multirow{4}{*}{5}
& $3$ & 19 & 5 & 0.6 & 19 & 7 & 0.64 & 20 & 8 & 0.66 \\
& $4$ & 17 & 5 & 3.88 & 18 & 6 & 4.13 & 18 & 8 & 4.16 \\
& $5$ & 17 & 5 & 26.1 & 17 & 5 & 23.6 & 19 & 6 & 28.3 \\
& $6$ & 17 & 5 & 154 & 16 & 6 & 149 & 16 & 6 & 154 \\
\hline
\end{tabular}}
\end{center}
\end{footnotesize}
\end{table}

From Tables \ref{Navier_Stokes_equations_table_Radau_Matlab_gamma1}--\ref{Navier_Stokes_equations_table_Radau_Matlab_gamma100}, we observe the robustness of the proposed preconditioner across a range of mesh sizes and viscosities: the linear solver converges in fewer than 27 iterations on average. The average number of iterations increases with the number of stages \( s \), though not drastically. Furthermore, we observe that as \( \gamma \) increases, the average number of FGMRES iterations decreases. The CPU times scale super-linearly with the system size, primarily due to the use of a direct solver for the augmented blocks. As we will demonstrate in the next section, employing an inexact solve for the diagonal blocks within the Gauss--Seidel method leads to a robust solver with CPU times scaling linearly with the system dimension.
Regarding Newton's method, the number of Newton iterations remains approximately constant, with a slight increase for small viscosities and a decrease as the mesh is refined.

We conclude this section by analyzing the dependence of the preconditioner on the viscosity $\nu$. To this end, we run our code with a refinement level of $l = 6$. Given the observed robustness of the preconditioner with respect to $\gamma$ in Tables \ref{Navier_Stokes_equations_table_Radau_Matlab_gamma1}--\ref{Navier_Stokes_equations_table_Radau_Matlab_gamma100}, we fix $\gamma = 1$ and vary the viscosity $\nu$. The results are reported in Table \ref{Navier_Stokes_equations_table_dependence_nu}. As shown, the number of linear iterations remains nearly constant, even in strongly convection-dominated regimes, while the number of Newton iterations increases slightly as the viscosity approaches zero.

\begin{table}[!ht]
\caption{Average FGMRES iterations ($\texttt{it}$), average Newton iterations ($\mathrm{Nit}$) per time step, and average CPU times per linear iteration for the block augmented Lagrangian preconditioner with $l=6$ and $\gamma=1$, for $\nu \in \left\{ \frac{1}{100}, \frac{1}{250}, \frac{1}{500}, \frac{1}{2500} \right\}$, for Radau IIA methods.}\label{Navier_Stokes_equations_table_dependence_nu}
\begin{footnotesize}
\begin{center}
\renewcommand{\arraystretch}{1.2}
{\begin{tabular}{|c||c|c|c||c|c|c||c|c|c||c|c|c|}
\hline
\multicolumn{1}{|c||}{} & \multicolumn{3}{c||}{$\nu=\frac{1}{100}$} & \multicolumn{3}{c||}{$\nu=\frac{1}{250}$}  & \multicolumn{3}{c||}{$\nu=\frac{1}{500}$} & \multicolumn{3}{c|}{$\nu=\frac{1}{2500}$}\\
\cline{2-13}
$s$ & $\texttt{it}$ & $\mathrm{Nit}$ & CPU & $\texttt{it}$ & $\mathrm{Nit}$ & CPU & $\texttt{it}$ & $\mathrm{Nit}$ & CPU & $\texttt{it}$ & $\mathrm{Nit}$ & CPU \\
\hline
\hline
3 & 16 & 3 & 83.8 & 16 & 4 & 86.3 & 16 & 4 & 91.6 & 19 & 8 & 111 \\
\hline
5 & 25 & 4 & 232 & 25 & 4 & 232 & 25 & 4 & 247 & 32 & 8 & 322 \\
\hline
\end{tabular}}
\end{center}
\end{footnotesize}
\end{table}

\subsection{Experiments with inexact block-diagonal solvers}
In this section, we study the effect of applying an inexact solver to the augmented blocks. For the spatial discretization, we use the Firedrake package \cite{Rathgeber_Ham_Mitchell_Lange_Luporini_Mcrae_Bercea_Markall_Kelley}, along with the \texttt{Irksome} package \cite{Farrell_Kirby_MarchenaMenendez, Kirby_MacLachlan} for Runge--Kutta time stepping. The recent extension of \texttt{Irksome} in \cite{Kirby_MacLachlan} enables a more general treatment of boundary conditions for DAEs, including support for discontinuous initial data (see \cite[Section 2.1]{Kirby_MacLachlan}). In the tests reported here, we use ODE-type boundary conditions, as DAE-type conditions may introduce a discontinuity at the start of time integration.

All tests are performed with $\gamma = 1$. For the finite elements, we use the exactly incompressible Scott--Vogelius $[P_4]^2$--$P_3^{\rm disc}$ pair \cite{Scott_Vogelius}.

FGMRES is used as the linear solver, with a block Gauss--Seidel preconditioner for the $(1,1)$-block \( \Phi + \gamma \Psi_1 \mathcal{W}^{-1} \Psi_2 \). Each diagonal block is approximated by one multigrid cycle following \cite{Farrell_Mitchell_Wechsung}, and the Schur complement is approximated as in \eqref{Schur_approx}. In the multigrid, we apply 10 FGMRES iterations preconditioned by an additive vertex-star relaxation on each fine level, and solve the coarse problem with LU decomposition.

All experiments are run on an Intel Core Ultra 7 165U processor with 64~GB RAM under Ubuntu 22.04.5 LTS.

\subsubsection{Lid-driven cavity}
First, we reproduce the tests for the time-dependent lid-driven cavity described in Section~\ref{sec_5_1_2}, using an inexact solver for the augmented diagonal blocks. Here, we examine the dependence of the solver on the family of Runge--Kutta methods employed. To this end, we vary the level of refinement $l$, the viscosity $\nu$, the number of stages $s$, and the Runge--Kutta family. We fix the number of time steps to $n_t = 50$.

For this test, the nonlinear iteration is performed until either the relative or absolute nonlinear residual is reduced by a factor of $10^{-5}$. FGMRES is run until the relative residual is reduced by $10^{-6}$, or the absolute residual by $10^{-5}$.

In Tables~\ref{Navier_Stokes_equations_table_Gauss}--\ref{Navier_Stokes_equations_table_Radau}, we report the average number of FGMRES iterations per nonlinear iteration, and the average CPU time (in seconds) per time step, for $s$-stage Gauss, Lobatto IIIC, and Radau IIA methods, with $s \in \{2, 3, 4\}$. Note that the reported CPU times include the evaluation of the nonlinear residual, discretization of the operators, and construction of the geometric multigrid on the augmented velocity blocks, in addition to the time required by the linear solvers. Table~\ref{Navier_Stokes_equations_table_Gauss} also reports the dimension (DoF) of each system solved.

\begin{table}[!ht]
\caption{Average FGMRES iterations ($\texttt{it}$) per non-linear iteration, average CPU times per time step, and degrees of freedon (DoF) with the block augmented Lagrangian preconditioner, for $\nu=\frac{1}{100}$, $\frac{1}{250}$, $\frac{1}{500}$, and $\frac{1}{2500}$ and a range of $l$, for $s$-stage Gauss, $s=2$, $3$, and $4$, and $\gamma=1$.}\label{Navier_Stokes_equations_table_Gauss}
\begin{footnotesize}
\begin{center}
\renewcommand{\arraystretch}{1.2}
{\begin{tabular}{|c|c||c||c|c||c|c||c|c||c|c|}
\hline
\multicolumn{1}{|c|}{} & \multicolumn{1}{c||}{} & \multicolumn{1}{c||}{} & \multicolumn{2}{c||}{$\nu=\frac{1}{100}$} & \multicolumn{2}{c||}{$\nu=\frac{1}{250}$} & \multicolumn{2}{c||}{$\nu=\frac{1}{500}$} & \multicolumn{2}{c|}{$\nu=\frac{1}{2500}$}\\
\cline{4-11}
$s$ & $l$ & DoF & $\texttt{it}$ & CPU & $\texttt{it}$ & CPU & $\texttt{it}$ & CPU & $\texttt{it}$ & CPU \\
\hline
\hline
\multirow{5}{*}{2} & $2$ & 6916 & 10 & 0.73 & 10 & 0.52 & 10 & 0.53 & 10 & 0.51 \\
 & $3$ & 27,140 & 10 & 1.85 & 10 & 1.79 & 10 & 1.82 & 9 & 1.72 \\
 & $4$ & 107,524 & 10 & 7.01 & 9 & 6.73 & 10 & 7.03 & 10 & 7.07 \\
 & $5$ & 428,036 & 10 & 27.4 & 10 & 27.5 & 10 & 27.4 & 9 & 28.8 \\
 & $6$ & 1,708,036 & 10 & 110 & 10 & 113 & 10 & 114 & 9 & 112 \\
\hline
\multirow{5}{*}{3} & $2$ & 10,374 & 27 & 2.24 & 31 & 1.98 & 28 & 1.82 & 30 & 1.93 \\
 & $3$ & 40,710 & 28 & 6.72 & 23 & 5.69 & 30 & 7.13 & 29 & 6.91 \\
 & $4$ & 161,286 & 25 & 24.4 & 29 & 28.1 & 31 & 29.9 & 24 & 26.6 \\
 & $5$ & 642,054 & 26 & 100.8 & 25 & 99.1 & 25 & 98.1 & 21 & 108 \\
 & $6$ & 2,562,054 & 24 & 379 & 25 & 392 & 28 & 431 & 20 & 417 \\
\hline
\multirow{5}{*}{4} & $2$ & 13,832 & 19 & 2.34 & 19 & 1.78 & 17 & 1.63 & 16 & 1.61 \\
 & $3$ & 54,280 & 14 & 5.01 & 18 & 6.36 & 18 & 6.18 & 17 & 6.18 \\
 & $4$ & 215,048 & 18 & 24.6 & 18 & 25.1 & 17 & 24.3 & 16 & 25.1 \\
 & $5$ & 856,072 & 18 & 100 & 18 & 100 & 18 & 99.1 & 14 & 104 \\
 & $6$ & 3,416,072 & 17 & 398 & 18 & 404 & 18 & 404 & 14 & 434 \\
\hline
\end{tabular}}
\end{center}
\end{footnotesize}
\end{table}

\begin{table}[!ht]
\caption{Average FGMRES iterations ($\texttt{it}$) per non-linear iteration and average CPU times per time step with the block augmented Lagrangian preconditioner, for $\nu=\frac{1}{100}$, $\frac{1}{250}$, $\frac{1}{500}$, and $\frac{1}{2500}$ and a range of $l$, for $s$-stage Lobatto IIIC, $s=2$, $3$, and $4$, and $\gamma=1$.}\label{Navier_Stokes_equations_table_Lobatto}
\begin{footnotesize}
\begin{center}
\renewcommand{\arraystretch}{1.2}
{\begin{tabular}{|c|c||c|c||c|c||c|c||c|c|}
\hline
\multicolumn{1}{|c|}{} & \multicolumn{1}{c||}{} & \multicolumn{2}{c||}{$\nu=\frac{1}{100}$} & \multicolumn{2}{c||}{$\nu=\frac{1}{250}$} & \multicolumn{2}{c||}{$\nu=\frac{1}{500}$} & \multicolumn{2}{c|}{$\nu=\frac{1}{2500}$}\\
\cline{3-10}
$s$ & $l$ & $\texttt{it}$ & CPU & $\texttt{it}$ & CPU & $\texttt{it}$ & CPU & $\texttt{it}$ & CPU \\
\hline
\hline
\multirow{5}{*}{2} & $2$ & 5 & 0.48 & 5 & 0.34 & 5 & 0.34 & 5 & 0.33 \\
 & $3$ & 5 & 0.97 & 5 & 1.11 & 5 & 1.07 & 4 & 0.98 \\
 & $4$ & 5 & 3.99 & 5 & 4.16 & 5 & 4.11 & 4 & 4.00 \\
 & $5$ & 5 & 16.3 & 5 & 15.5 & 4 & 14.9 & 4 & 17.1 \\
 & $6$ & 4 & 60.4 & 4 & 59.8 & 4 & 60.8 & 4 & 59.9 \\
\hline
\multirow{5}{*}{3} & $2$ & 7 & 0.89 & 7 & 0.62 & 7 & 0.60 & 7 & 0.61 \\
 & $3$ & 8 & 2.21 & 7 & 2.22 & 7 & 2.13 & 6 & 2.00 \\
 & $4$ & 7 & 8.48 & 7 & 8.19 & 7 & 8.23 & 5 & 7.57 \\
 & $5$ & 7 & 33.3 & 6 & 31.7 & 5 & 28.3 & 5 & 33.2 \\
 & $6$ & 7 & 131 & 7 & 132 & 6 & 122 & 5 & 134 \\
\hline
\multirow{5}{*}{4} & $2$ & 9 & 1.45 & 9 & 1.01 & 8 & 1.01 & 9 & 1.01 \\
 & $3$ & 8 & 3.41 & 9 & 3.57 & 8 & 3.42 & 8 & 3.22 \\
 & $4$ & 8 & 13.1 & 8 & 13.3 & 8 & 13.2 & 7 & 13.1 \\
 & $5$ & 8 & 53.3 & 8 & 52.3 & 8 & 51.4 & 6 & 55.3 \\
 & $6$ & 8 & 212 & 8 & 212 & 7 & 205 & 6 & 232 \\
\hline
\end{tabular}}
\end{center}
\end{footnotesize}
\end{table}

\begin{table}[!ht]
\caption{Average FGMRES iterations ($\texttt{it}$) per non-linear iteration and average CPU times per time step with the block augmented Lagrangian preconditioner, for $\nu=\frac{1}{100}$, $\frac{1}{250}$, $\frac{1}{500}$, and $\frac{1}{2500}$ and a range of $l$, for $s$-stage Radau IIA, $s=2$, $3$, and $4$, and $\gamma=1$.}\label{Navier_Stokes_equations_table_Radau}
\begin{footnotesize}
\begin{center}
\renewcommand{\arraystretch}{1.2}
{\begin{tabular}{|c|c||c|c||c|c||c|c||c|c|}
\hline
\multicolumn{1}{|c|}{} & \multicolumn{1}{c||}{} & \multicolumn{2}{c||}{$\nu=\frac{1}{100}$} & \multicolumn{2}{c||}{$\nu=\frac{1}{250}$} & \multicolumn{2}{c||}{$\nu=\frac{1}{500}$} & \multicolumn{2}{c|}{$\nu=\frac{1}{2500}$} \\
\cline{3-10}
$s$ & $l$ & $\texttt{it}$ & CPU & $\texttt{it}$ & CPU & $\texttt{it}$ & CPU & $\texttt{it}$ & CPU \\
\hline
\hline
\multirow{5}{*}{2} & $2$ & 3 & 0.41 & 3 & 0.23 & 3 & 0.26 & 3 & 0.25 \\
 & $3$ & 2 & 0.64 & 3 & 0.72 & 2 & 0.67 & 2 & 0.64 \\
 & $4$ & 2 & 2.45 & 3 & 3.15 & 2 & 2.45 & 2 & 2.43 \\
 & $5$ & 2 & 9.68 & 2 & 9.39 & 2 & 9.11 & 2 & 10.7 \\
 & $6$ & 2 & 37.4 & 2 & 36.5 & 2 & 36.3 & 2 & 39.4\\
\hline
\multirow{5}{*}{3} & $2$ & 6 & 0.85 & 7 & 0.59 & 6 & 0.59 & 6 & 0.58 \\
 & $3$ & 6 & 2.01 & 7 & 2.05 & 6 & 1.95 & 6 & 1.92 \\
 & $4$ & 6 & 7.93 & 6 & 7.61 & 6 & 7.74 & 5 & 7.74 \\
 & $5$ & 6 & 31.1 & 6 & 30.5 & 6 & 29.7 & 5 & 31.8 \\
 & $6$ & 6 & 124 & 6 & 124 & 6 & 119 & 5 & 125 \\
\hline
\multirow{5}{*}{4} & $2$ & 4 & 1.27 & 4 & 0.67 & 4 & 0.65 & 4 & 1.23 \\
 & $3$ & 4 & 2.15 & 4 & 2.15 & 4 & 2.11 & 4 & 2.01 \\
 & $4$ & 4 & 8.21 & 4 & 8.22 & 4 & 7.91 & 3 & 8.18 \\
 & $5$ & 4 & 32.7 & 4 & 32.1 & 3 & 31.1 & 3 & 35.8 \\
 & $6$ & 4 & 129 & 3 & 128 & 3 & 126 & 3 & 137 \\
\hline
\end{tabular}}
\end{center}
\end{footnotesize}
\end{table}

From Tables~\ref{Navier_Stokes_equations_table_Gauss}--\ref{Navier_Stokes_equations_table_Radau}, we observe the expected robustness of the preconditioner when applying a multigrid routine as an inexact solver to the augmented blocks. The outer linear solver converges within at most 31 iterations to the prescribed tolerance, even in convection-dominated regimes (see, for example, the case with $\nu = \frac{1}{2500}$). Moreover, the use of a multigrid routine yields CPU times that scale linearly with the system size. Finally, we note that Newton's method converges to the prescribed tolerance in an average of 1 to 2 iterations per time step.

\subsubsection{Flow in a Channel around a Cylinder}
As a final test, we present results for flow in a 2D channel around a cylinder \cite{Schäfer1996}. The domain is given by $\Omega = [0, 2.2] \times [0, 0.41] \setminus B_r(0.2, 0.2)$, where $B_r(0.2, 0.2)$ denotes a disk of radius $r = 0.05$ centered at $(0.2, 0.2)$. We impose no-slip boundary conditions on the cylinder and on the top and bottom walls, and prescribe the following Poiseuille inflow on the left boundary:
\begin{displaymath}
u(x_1, x_2, t) = \left[ \frac{2 x_2 (H - x_2)}{H^2}, 0 \right]^\top,
\end{displaymath}
where $H = 0.41$. Natural outflow boundary conditions are imposed on the right boundary. We set $T = 1$ and divide the time interval into $n_t = 100$ subintervals.

For the spatial discretization, we start from an unstructured grid with 972 triangular elements and perform three levels of refinement, resulting in a mesh with 217{,}392 degrees of freedom for the velocity space and 134{,}560 for the pressure space. We solve the problem for $\nu = \frac{1}{100}$, $\frac{1}{500}$, and $\frac{1}{1000}$.

The nonlinear iteration is run until a reduction of $10^{-4}$ is achieved in either the relative or absolute residual. FGMRES is run until the absolute residual is reduced by $10^{-4}$ or the relative residual by $10^{-6}$. In Table~\ref{Navier_Stokes_equations_table_Radau_cylinder}, we report the average number of FGMRES iterations per nonlinear iteration, the average number of Newton iterations per time step, and the average CPU time per time step, along with the dimensions of the systems solved, for $s$-stage Radau IIA methods with $s = 2$, 3, and 4.

\begin{table}[!ht]
\caption{Average FGMRES iterations ($\texttt{it}$), average Newton iterations ($\mathrm{Nit}$), average CPU times, and degrees of freedon (DoF) per time step with the block augmented Lagrangian preconditioner, with $\nu=\frac{1}{100}$, $\frac{1}{500}$, and $\frac{1}{1000}$, for $s$-stage Radau IIA, $s=2$, $3$, and $4$, and $\gamma=1$.}\label{Navier_Stokes_equations_table_Radau_cylinder}
\begin{footnotesize}
\begin{center}
\renewcommand{\arraystretch}{1.2}
{\begin{tabular}{|c||c||c|c|c||c|c|c||c|c|c|}
\hline
\multicolumn{1}{|c||}{} & \multicolumn{1}{c||}{} & \multicolumn{3}{c||}{$\nu=\frac{1}{100}$} & \multicolumn{3}{c||}{$\nu=\frac{1}{500}$} & \multicolumn{3}{c|}{$\nu=\frac{1}{1000}$}\\
\cline{3-11}
$s$ & DoF & $\texttt{it}$ & $\mathrm{Nit}$ & CPU & $\texttt{it}$ & $\mathrm{Nit}$ & CPU & $\texttt{it}$ & $\mathrm{Nit}$ & CPU \\
\hline
\hline
2 & 703,904 & 23 & 2 & 23.5 & 20 & 1 & 20.4 & 21 & 2 & 22.1 \\
\hline
3 & 1,055,856 & 28 & 1 & 34.3 & 25 & 1 & 32.6 & 25 & 1 & 32.4 \\
\hline
4 & 1,407,808 & 37 & 1 & 60.7 & 33 & 1 & 53.2 & 29 & 1 & 48.5 \\
\hline
\end{tabular}}
\end{center}
\end{footnotesize}
\end{table}

From Table~\ref{Navier_Stokes_equations_table_Radau_cylinder}, we again observe the robustness of the preconditioner when applied to flows in open domains. The solver converges to the prescribed tolerance in at most 37 iterations for the parameter values considered here.

\section{Conclusions}\label{sec_6}
In this work, we investigated the Runge--Kutta discretization of the nonstationary incompressible Navier--Stokes equations. At each time step, the stage equations of the Runge--Kutta method yield a nonlinear system, which we solve using Newton’s method. The resulting linear systems exhibit a generalized saddle-point structure. Leveraging classical saddle-point theory, we developed an augmented Lagrangian-based block triangular preconditioner tailored for these systems.
Through a series of numerical experiments, we demonstrated the robustness and efficiency of the proposed preconditioner with respect to key parameters, including viscosity, mesh size, time step, and the number of Runge--Kutta stages. Notably, the preconditioner remains effective even when the inverses of the augmented diagonal blocks are applied inexactly using a multigrid routine.

Future work will focus on a rigorous theoretical analysis of the proposed method and its extension to more complex flow systems.

\section*{Acknowledgements}
S.L. greatfully acknowledges Scott MacLachlan for useful discussions on the Firedrake system and the \texttt{Irksome} package. S.L. is a member of Gruppo Nazionale di Calcolo Scientifico (GNCS) of the Istituto Nazionale di Alta Matematica (INdAM). M.O. was supported in part by the U.S. National Science Foundation under awards DMS-2309197 and DMS-2408978.

%
%
%

\end{document}